\documentclass[12pt]{article}
\usepackage{amsfonts}
\begin{document}
\begin{center}
{\bf THE LIFE, WORK AND LEGACY OF P. L. CHEBYSHEV}
\end{center}
\begin{center}
{\bf N. H. BINGHAM}
\end{center}
\noindent {\bf Abstract} \\
\indent We survey briefly the life and work of P. L. Chebyshev, and his ongoing influence.  We discuss his contributions to probability, number theory and mechanics, his pupils and mathematical descendants, and his role as the founding father of Russian mathematics in general and of the Russian school of probability in particular. \\

{\bf 1. Life} \\

 Pafnuty L'vovich Chebyshev was born on 16 May 1821 in Okatovo in the Kaluga district of Russia.  His father was a Russian nobleman and wealthy landowner.  He was at first home-educated; the family moved to Moscow in 1832, to obtain better tuition.  He entered the University of Moscow in 1837, aged 16 (unusually early now, entry to university at such an age was not unusual in former times, for example in England's two ancient universities).  He graduated in 1843.  He wrote his Master's thesis in 1846, {\sl An essay on the elementary analysis of the theory of probability}.  While at Moscow, Chebyshev was taught by N. D. Brashman (1796-1866), whom he greatly respected.  Brashman's interests included mechanics, mechanical engineering and hydraulics.  This breadth of interests, and in particular the interest in mechanical engineering, may have influenced Chebyshev's interests in such directions (\S 6).  \\
\indent As he found no academic post in Moscow, Chebyshev moved to St. Petersburg, where he obtained his {\it venia legendi} (right to teach) and became a lecturer in the University.  He took his doctorate in 1849, on the theory of numbers.  He became an extraordinary professor of mathematics in 1850, and an ordinary (full) professor in 1860.  He became a junior academician at the St. Petersburg Academy of Sciences with a chair in applied mathematics in 1853, and an ordinary (again, full) academician in 1859.  He retired from teaching at the University in 1882, but continued to do research at the Academy all his life.  He received many foreign honours. \\
\indent Chebyshev suffered from childhood from a medical condition (now called Trendelenburg's gait, named after the first person to study it in 1895), which made one leg longer than the other.  He had a limp, and walked with a stick.  Of course, this denied him much by way of childhood games and sport, and also the army career his parents had intended for him (his father was an Army officer who fought Napoleon); instead, he turned to mathematics.  One wonders what effect modern physiotherapy might have achieved here if applied in good time. \\
\indent Chebyshev began his career publishing in Russian.  Russia had a scientific tradition: St. Petersburg University was founded in 1724 by Peter the Great, and the St. Petersburg Academy in 1725 by his widow Catherine I; Daniel Bernoulli was a professor at the Academy from 1724-33; Euler was at the Academy from 1727-41, succeeding Bernoulli.  But as Chebyshev became more established, he realised that if he wanted his work to receive proper recognition, he should publish in a Western European language.  The natural choice was French, the preferred language of the Russian nobility in tsarist times (as Greek was with the Roman aristocracy in ancient times).  Incidentally, Chebyshev's use of French gave rise, via the usual French transliteration from Cyrillic to Roman, to the symbol $T_n$ for the Chebyshev polynomials (of the first kind) -- the place where many students first meet his name -- from Tchebychev or Tchebycheff (the T being needed in French to harden the first consonant).  \\
\indent Chebyshev was regarded as a good teacher; his lectures were lively, stimulating and appreciated.  His pupil Lyapunov took careful notes of his lectures on probability for 1879-80; see [She94], [Kry36].  His two favourite courses were on number theory and probability theory; he taught each thirty-one times. \\  

{\bf 2. Number theory} \\

 In Edmund Landau's classic [Lan09] on the distribution of prime numbers, he begins his account of Chebyshev (I.4, `Tschebyschef', p.11-29): `Im allgemeinen Primzahlproblem hat nach Euklid erst Tschebyschef die ersten weiteren sicheren Schritte gemacht und wichtige S\"atze bewiesen' (In the general prime-number problem, after Euclid  Chebyshev was the first to make great and certain steps and to prove important results). \\  
\indent We begin with Bertrand's postulate of 1845 [Ber45]: for any positive integer $n \geq 1$ the interval $[n, 2n]$ contains at least one prime (for a very short proof, see Ramanujan [Ram19]).  \\
\indent Chebyshev proved Bertrand's postulate in 1852 [Che52], by deducing it from the following effective (quantitative) estimate for the prime-counting function
$$
\pi (x) := \sum_{p \leq x} 1
$$
for the primes $p \leq x$:
$$
\{ c_1 + o(1) \} x/\log x
\leq \pi(x) \leq
\{c_2 + o(1) \} x/\log x
$$
with
$$
c_1 = \log (2^{1/2} 3^{1/3} 5^{1/5}/{30}^{1/30}) = 0.92129 \cdots, \quad
c_2 = 6 c_1/5 = 1.10555 \cdots .
$$
In particular, {\it if} the limit $ \pi(x)/(x/\log x)$ exists, it is $1$:
$$
                 \pi(x) \sim x/\log x \qquad (x \to \infty).                  \leqno(PNT)
$$ 
This would prove the Prime Number Theorem (PNT), conjectured in 1792 by Gauss (then aged 15) on numerical grounds.  See e.g. Tenenbaum [Ten15, \S 1.2].  Otherwise put, {\it if} $\pi(x) \sim c x/\log x$ for some constant $c$, then $c = 1$:
$$
\liminf \pi(x)/(x/\log x) \leq 1 \leq \limsup \pi(x)/(x/ \log x)
$$
[Ten15, \S 1.7].  For another treatment here, see Rose [Ros88 \S 12.1]. \\
\indent Write
$$
\Lambda(n) := \log p \quad (n = p^k, k \geq 1), \quad 0 \quad (n \neq p^k)
$$
for the von Mangoldt function $\Lambda$.  As well as $\pi(x)$ above, Chebyshev gave us two other summatory functions,
$$
\psi(x) := \sum_{n \leq x} \Lambda(n), \qquad
\theta(x) := \sum_{p \leq x} \log p
$$
[Ten, \S 2.6].  These are ubiquitous in prime-number theory.  For example,
$$
\psi(x) \sim x
$$
is `elementarily equivalent to' (can easily be proved equivalent to) $(PNT)$ ([Lan09]; [Ten15, \S 3.6]). \\

{\bf 3. Probability theory} \\

   We turn now to the field where Chebyshev's name is best known -- probability theory.  His pre-eminence here is well exemplified by the careful twenty-page summary of his work in Maistrov's book [Mai74, IV.1], and by the saying of an eminent probabilist of our time: that over all time, the three key names are Bernoulli, Chebyshev, Kolmogorov.  \\

1. {\it Weak laws of large numbers} \\
\indent The early paper [Che46] of Chebyshev's on probability, following his Master's thesis, also of 1846 [Mai74, 191-195], is on weak laws of large numbers for Bernoulli trials with varying success probabilities $p_i$, following work by Poisson in his book of 1837 (the origin of the Poisson distribution).  It is considered in some detail in Stigler [Stig86, 182-186] and Maistrov [Mai74, 195-198].  Maistrov points out that independence, which Chebyshev assumes, is never stated explicitly. \\ 
 
 2. {\it The Bienaym\'e-Chebyshev inequality} \\
\indent The main  result of the paper [Che67] of 1867 -- that for a random variable with two moments finite
$$
\mathbb{P} (|X - \mathbb{E}[X]| \geq t) \leq var (X)/t^2 \quad (t > 0)
$$
-- is the single result in probability for which Chebyshev's name is best known, and is ubiquitous in textbooks as `Chebyshev's inequality' (e.g. [Fel, IX.6], which is where the author met it).   \\
\indent The result was in fact proved fourteen years earlier by I. J. Bienaym\'e [Bie53].  Remarkably, when Chebyshev re-discovered the result, and it was to appear in Liouville's journal {\sl Journal de Math\'ematiques Pures et Appliqu\'es}, Liouville arranged for Bienaym\'e's paper to be printed immediately before Chebyshev's.  It is even more remarkable that despite this the name `Chebyshev's inequality' stuck. \\
\indent Some authors use the name Bienaym\'e's inequality.  See for example Jean-Pierre Kahane's remarkable book on random series, [Kah85].  This book, incidentally, can be read as an extended essay on the power of two results, this inequality and Fubini's theorem, used systematically together. \\ 
\indent Weak laws of large numbers such as those above were as far as laws of large numbers could develop in the nineteenth century: even the formulation of strong laws had to await measure theory and the twentieth century. \\

3. {\it The method of moments} \\
\indent Even when all the moments of a probability distribution exist, they may or may not determine the distribution uniquely (the question of when they do is the {\it moment problem}, \S 5).  When all moments exist, convergence of moments implies convergence of distributions (see e.g. Billingsley [Bil95, \S 30]).  This is the {\it method of moments} in probability theory.  It was used by Chebyshev [Che74], motivated by the central limit theorem (below).  It is linked to his work on continued fractions (\S 4.4), Gaussian quadrature (\S 4.7), the Chebyshev separation theorem, and Chebyshev systems (\S 4.8).  For commentary, see Krein [Kre51], Kjeldsen [Kje93]. \\

4. {\it The central limit theorem} \\
\indent The next category of limit theorems in probability theory is that of central limit theorems.  The first version here, for Bernoulli trials, goes back to de Moivre in 1733, in the paper {\sl Approximatio ad Summam Terminorum Binomii} ${\overline{a + b}}^n$ {\sl in Seriem expansi}, circulated but not published, and incorporated in the second, 1738, edition of his book {\sl The Doctrine of Chances}.  Here the standard normal (or Gaussian) law with density
$$
\phi(x) := \frac{1}{\sqrt{2 \pi}} \exp \{ - \frac{1}{2} x^2 \}
$$
enters mathematics and science forever. (Note that Stirling's formula, intimately connected with this, was proved in 1730.)  See [Stig86, 78-88] for a full treatment.  Wide generalization and application of de Moivre's work, most notably to astronomy and celestial mechanics, had to wait for Laplace, over the period 1777-1812 [Stig86, 117-138], Legendre, 1805-1820, and Gauss, 1795-1809.  See Stigler [Stig86, Ch. 4] for the resulting `Gauss-Laplace synthesis', and Seal [Sea67].  The name `de Moivre-Laplace limit theorem' is usually given to the Bernoulli case of the central limit theorem.  Usage varies regarding the  name given to the density above: Laplacian is common in France, Gaussian in most countries; the term normal is due to Karl Pearson around 1900, and is used synonymously with Gaussian in English. \\
\indent In his {\sl Acta Mathematica} paper [Che91] of 1891 (Russian, 1887), Chebyshev treats both the weak law of large numbers and the central limit theorem, the latter using the {\it method of moments}.  Despite some technical deficiencies in the paper, this was a valuable theoretical advance; so too were the error estimates Chebyshev obtained.  Maistrov [Mai74, 206, 207] quotes comments by Kolmogorov and by Gnedenko on this paper, both written in 1948, and both pre-figuring the focus of their own forthcoming book, the classic [GneK] of 1949 (English translation 1954), which gives the definitive solution of the general central-limit problem.  Kolmogorov continues [Mai74, 207]: ``P. L. Chebyshev impelled Russian probability into first place in the world". \\
\indent This assessment is echoed by that of Khinchin [Mai74, 208]: that from the second half of the nineteenth century, Russia was ``the only country in which the mathematical foundations of probability theory were cultivated with the seriousness it deserved in view of its prominent role in the natural sciences and engineering.  It is entirely due to the works of Chebyshev that the Russian school of probability attained this exceptional position".  \\
\indent One can only agree. \\ 

{\bf 4. Analysis} \\

{\it 1. Polynomials of best approximation} \\
\indent Chebyshev's work here [Che58] is his best-known in analysis; for a modern textbook account, see e.g. Natanson [Nat64, Ch. II.2].  The main result is {\it Chebyshev's Alternation Theorem}: for a function $f \in C[a,b]$, a polynomial $p_n$ of degree $n$ is a polynomial of {\it best approximation} (uniformly -- in the uniform norm) if and only if there are $n+2$ points $x_i$, $a \leq x_0 < \cdots < x_{n+1} \leq b$, where $\pm (-1)^i [f(x_i) - p_n(x_i)]$ takes its maximum value $\Vert f - p_n \Vert$, {\it alternately}.  In Natanson's terminology, the $(+)$-points and the $(-)$-points {\it alternate}.  Such a set of points is called a {\it Chebyshev alternant}.\\
\indent From this same paper emerged Chebyshev's work on orthogonal polynomials in general, the orthogonal polynomials $T_n$ that bear his name, and his work on continued fractions; see below. \\
\indent One may mention here the other main result on polynomial approximation, the Weierstrass approximation theorem (1885), and Bernstein's probabilistic proof via his polynomials (1912/13) [Nat64, I.1], [Lor53, I.1.1]. \\
 
{\it 2. Best rational approximation} \\
\indent The Alternation Theorem is extended to {\it rational functions} in the long sequel [Che59].  For a textbook account, see Akhiezer [Akh56, II.34].\\
    
{\it 3. Explicit integration} \\
\indent It is all too easy to write down a simple integrand, and (as with the standard normal density above) find that the only exact expression for its integral is `integrand plus integral sign'.  Chebyshev applied himself to such problems of `integration in finite terms' [ButJ99].  In his thesis he proved a conjecture of Abel from 1826: if
$$
\int (\rho(x)/R(x)) dx \qquad(\rho, R \  \hbox{polynomials})
$$
is expressible by logarithms, then
$$
I := \int \frac{\rho(x)}{\sqrt{R(x)}} dx = c \log \frac{p + q \sqrt{R}}{p - q \sqrt{R}}
$$
with $p$, $q$ entire functions and $c$ constant. \\
\indent The corresponding problem with $m$th roots was considered by Abel and Liouville, who found that if $I$ is expressible in finite form, then
$$
I = U + c_0 \log V_0 + \cdots + c_n \log V_n,
$$
with $U, V_0, \cdots, V_n$ rational functions of $x + R(x)^{1/m}$.  In 1852 Chebyshev found $U$, and how many $c_i \log V_i$ terms are needed. \\
\indent Hardy began his career as a noted expert on the integral calculus; his series {\sl Notes on some points in the integral calculus} ran to 49 papers (1901-1929; numbered in Roman numerals; `N.I.C.' in his Collected Works).  Hardy's first book was an early Cambridge Tract on this subject [Har05].  It cites five papers by `P. Tschebyschef', in the period 1853-1861. \\
\indent This area was the subject of a later monograph by Ritt [Rit48]. \\
 
{\it 4. Continued fractions} \\
\indent Continued fractions can be traced back to the Euclidean algorithm, to the manual extraction of square roots (still taught when my father, b. 1904, was at school), Diophantine equations etc.  For the early theory, see e.g. Brezinski [Bre91].  Apart from the arbitrariness of the base 10 of decimals, the criterion for rationality is simpler for a continued-fraction expansion (finiteness, i.e. termination) than for a decimal (termination {\it or} recurrence).  The criterion for recurrence of a continued-fraction expansion is (Lagrange's theorem) being a {\it quadratic irrational} (see e.g. [HarW, Ch. X]), a result of real number-theoretic importance.  Such results suggest that continued fractions give the natural way to expand real numbers.  Indeed, as ${\mathbb{N}}^{\mathbb{N}}$ is in bijection with the irrationals under the continued-fraction map, it is used in descriptive set theory to {\it encode} the irrationals; see e.g. Sierpinski [Sie28, \S 35], [Sie64, VIII.2]; [RogJ,  \S 2.1].  \\
\indent An important and very early result was Brouncker's continued-fraction for $\pi$ (actually $4/\pi$), of 1655:
$$
4/\pi = 1 + \frac{1^2}{2 \ +} \ \frac{3^2}{2 \ +} \ 
\frac{5^2}{2 \ +} \ \frac{7^2}{2 \ +} \ \cdots 
$$
(where the $+$ at the end of each denominator means `start the new `fraction within a fraction' here' [Bre91, 1.4]).  This was incorporated into Wallis's {\sl Arithmetica Infinitorum} of 1656 (it is related to Wallis's infinite-product expansion for $\pi$).  Euler in turn was deeply influenced by Wallis's book, and obtained (1737) the corresponding continued fraction for $e$ [Bre91, 4.1]:
$$
e = 2 + \frac{1}{1 \ +} \ \frac{1}{2 \ +} \ \frac{1}{1 \ +} \ \frac{1}{1 \ +} \ \frac{1}{4 \ +} \ \frac{1}{1 \ +} \ \frac{1}{1 \ +} \ \frac{1}{6 \ +}  \cdots 
$$

{\it 5. Orthogonal polynomials and continued fractions} \\
\indent The theory of orthogonal polynomials $p_n$ (with respect to an orthogonality measure $\mu$) may be founded on the {\it three}-term recurrence formula with {\it two} sequences of coefficients, $(a_k)$ and $(b_k)$ (see e.g. Akhiezer [Akh65, p.4]):
$$
x P_k(x) = b_{k-1} P_{k-1}(x) + a_k P_k(x) + b_k P_{k+1}(x) \quad (k = 0,1,\cdots; \ b_{-1} = 0).
$$
Given the measure $\mu$ there exist polynomials $p_n$ orthogonal to it and satisfying this recurrence (see e.g. Szeg\H{o} [Sze59, Th. 3.2.1]), and conversely (Favard's theorem: [Sze59, 43]). \\
\indent Call the continued fraction [Akh65, (1.36)]
$$
\frac{1}{x - a_0 \ -} \ \ \frac{b_0^2}{x - a_1 \ -} \ \ 
\frac{b_1^2}{x - a_2 \ -} \   \cdots
$$
that {\it corresponding to} the recurrence relation above.  Write its convergents (partial products) as $R_n(x)/S_n(x)$.  Then [Sze59, p.55]
$S_n(x) = \sqrt{c_0} p_n(x)$: {\it the orthogonal polynomials are the denominators of the continued fraction}.  This is the essence of the link between the two. \\
\indent Accompanying this continued fraction is an infinite {\it Jacobi matrix}, a tridiagonal matrix
$$
\cal{J} :=
\left(
\begin{array}{cccccc}
a_0 & b_0 & 0   & 0   & 0 & \ \\
b_0 & a_1 & b_1 & 0   & 0 & \ \\
0   & b_1 & a_2 & b_2 & 0 & \ \\
0   & 0   & \ddots & \ddots & \ddots & \
\end{array}
\right) 
$$
[Akh65, p.2]. \\
\indent Szeg\H{o} [Sze59, 3.5] cites eight papers of Chebyshev here, beginning with [Che58], from 1858 to 1885; Khrushchev [Khr08] cites six, from 1854 to 1857; Brezinski [Bre91] cites 34, from 1855 to 1907. \\

{\it 6. Orthogonal polynomials} \\  
\indent Though continued fractions go back much earlier (above), and special families of orthogonal polynomials go back to the 18th century (Legendre polynomials, 1785), the general theory of orthogonal polynomials stems from Chebyshev in 1858 [Che58] (followed by further special cases: Hermite, 1864, Laguerre, 1879, Gegenbauer, 1884, etc.) \\ 
\indent The polynomials for which Chebyshev is remembered are
$$
T_n := \cos (n \ {\cos}^{-1}x): \quad T_n(cos \ \theta) = \cos n \theta \quad (n = 0,1,2,\cdots)
$$
One can easily check by de Moivre's theorem that $T_n$ is indeed a polynomial of degree $n$.  By the orthogonality property of the cosine function, the $T_n$ are easily seen to be orthogonal on $[-1,1]$ under the measure $dx/\sqrt{1 - x^2}$.  From the alternation of the maxima and minima of the cosine, these extrema of $T_n$ form a Chebyshev alternant.  So by Chebyshev's Alternation Theorem, $T_n$ gives the polynomial of degree $n$ of least absolute deviation from zero on $[-1,1]$.  For a full account of the importance of the $T_n$ in uniform approximation, see again [Nat64], Davis [Dav63, \S 3.3], or Lorentz [Lor86, Ch. 2] (for approximation in mean, see [Nat65a]).  It is worth noting that [Che58] contains the Christoffel-Darboux identity ([Sze59, \S 3.2]; Christoffel 1858, Darboux 1878). \\

{\it 7. Gaussian quadrature} \\
\indent With the points $x_i \in (a,b)$ specified in advance, the classical Newton-Cotes (or Lagrange) quadrature formula
$$
I(f) := \int_a^b f(x) d\mu(x) = \sum_1^n {\lambda}_k f(x_k),
$$
based on the Newton (divided-difference) or Lagrange interpolation formulae,  is exact for all polynomials of degree $n$ whenever the weights ${\l}_k$ are suitably chosen: that is, whenever it is exact for the Lagrange interpolation polynomials ${\ell}_1, \cdots, {\ell}_n$.  But if we choose the points $x_k$ suitably, such a formula can be made exact for all polynomials of degree $2n-1$.  For this, we choose the $x_k$ to be the $n$ roots of the $n$th orthogonal polynomial $P_n$ with respect to the measure $\mu$ on $[a,b]$; see e.g. [Sze59, \S 3.4].  This gives the {\it Gaussian quadrature} (or Gauss-Jacobi mechanical quadrature) for $\mu$ on $[a,b]$, [Nat65b].  The weights ${\lambda}_k$ here are called the {\it Christoffel numbers}.  They sum to the mass $\mu$ gives to the interval:
$$
\sum_1^n {\lambda}_k = \mu(b) - \mu(a).
$$
The zeros (written $x_1 < \cdots < x_n$ for convenience) of each such $P_n$ are real, distinct, and lie in $(a,b)$.  The zeros of $P_n$ and $P_{n+1}$ interlace [Sze59, \S 3.3]. \\
\indent Imagine that one starts at the left end-point $a =: y_0$, and moves to the right, watching the mass of the measure $\mu$ accumulate.  The point where it reaches ${\lambda}_1$ is $y_1$; the point where further mass ${\lambda}_2$ has accumulated is $y_2$, etc., and finally $y_n = b$.  The {\it Chebyshev separation theorem} [Che74] (stated by Chebyshev in 1874, proved later by A. A. Markov and Stieltjes) states that the points $x_1 < \cdots < x_n$ and $y_1 < \cdots <  y_{n-1}$ interlace.  Then the interpolation formula becomes
$$
\int_a^b f(x) d\mu(x) \sim \sum_1^n f(x_k) [\mu(x_k) - \mu(x_{k-1}],
$$
a Riemann-Stieltjes sum for the integral.  Thus the sum will converge to the integral as the degree $n$ increases, under suitable conditions, giving {\it convergence of quadrature processes}.  For $d\mu(x) = w(x) dx$, this holds for $f \in C[a,b]$, and more generally (Stieltjes; see e.g. [Dav63, Cor. 14.4.7]). \\
\indent It is a mark of the importance of the Chebyshev separation theorem that Szeg\H{o} gives three separate proofs of it [Sze59, 3.41].  As Szeg\H{o} points out, there is a degree of non-uniqueness in the $y_i$ as introduced above because of atoms in the measure $\mu$.  As he also points out, this does not affect the quadrature formula, but for more detail here see [Akh65, Th. 2.5.4]. \\
\indent For links between such quadrature formulae (extended to `quasi-orthogonal polynomials') and continued fractions, see [Akh65, \S 1.4]. \\ 

{\it 8. Chebyshev systems} \\
\indent The title `On the limiting values of integrals' of [Che74] (the term is due to him: [Akh65, vi]) would perhaps be more helpfully rendered nowadays along the lines of `On the bounds of integrals subject to constraints'. \\
\indent A {\it Chebyshev system} on $[a,b]$ is a set $u_0, \cdots, u_n$ of continuous real-valued functions for which all determinants
$$
\left| \begin{array}{cccc}
u_0(t_0) & u_0(t_1) & \cdots & u_0(t_n) \\ 
u_1(t_0) & u_1(t_1) & \cdots & u_1(t_n) \\ 
\vdots   & \vdots   & \      & \vdots   \\
u_n(t_0) & u_n(t_1) & \cdots & u_n(t_n)
\end{array} \right|
$$
are strictly positive whenever $a \leq t_0 < \cdots < t_n \leq b$.  The classic example is the set of powers $u_k(t) = t^k$, when the determinant reduces to the Vandermonde determinant with value
$$
\prod_{i < j} (t_j - t_i).
$$
\indent M. G. Krein ([Kre51]; Boas [Boa52]) gives a detailed account of `the ideas of P. L. \v Ceby\v sev and A. A. Markov in the theory of limiting values of integrals and their further developments', beginning with a detailed historical account.  We are to estimate 
$$
\int_{\xi}^{\eta} \Omega d \sigma 
$$ 
($a \leq \xi < \eta \leq b$, $\Omega$ continuous, $\sigma$ non-decreasing), with the first $n+1$ `generalized moments' prescribed:
$$
\int_a^b u_k d \sigma = c_k \quad (k = 0, 1, \cdots, n).
$$
The classical case of the powers is noted above.  There one can use continued fractions; here new methods are needed.  Krein uses the `method of maximal mass'.  For a textbook account of this, the Markov-Krein theorem and much else, we refer to Karlin's monograph [Kar66].     \\

{\bf 5. The moment problem} \\

 Chebyshev's work on continued fractions and orthogonal polynomials leads naturally on to the moment problem.  The main developments lead us beyond Chebyshev's lifetime, so this concerns his legacy, although it builds on his work.  The two classics here are Akhiezer [Akh65] and Shohat and Tamarkin [ShoT]; we use mainly the first (admirably clear, though longer). \\
\indent The {\it Hamburger moment problem} applies to the line: given a sequence of reals $(s_k)_0^{\infty}$ ($s_0 = 1$), find a measure $\mu$ with the $s_k$ as moments:
$$
s_k = \int_{-\infty}^{\infty} u^k d\mu(u) \quad (k = 0,1,\cdots)   \leqno(Ham)
$$
(thus $\int d\mu = s_0 = 1$: $\mu$ is a probability measure).  To avoid trivial cases, we restrict to $\mu$ with infinitely many points of increase.    Then (Hamburger, 1920, 1921; [Akh65, Th. 2.1.1]) such a $\mu$ exists if and only if $(s_k)$ is {\it positive}, i.e., all the Hankel quadratic forms
$$
\sum_{i,j = 0}^n s_{i+j} x_i x_j \quad (n = 0,1,2,\cdots)
$$
are positive definite, that is, that all the determinants
$$
D_k := 
\left| \begin{array}{cccc}
s_0 & s_1 & \cdots & s_k \\ 
s_1 & s_2 & \cdots & s_{k+1} \\                  
\vdots   & \vdots   & \      & \vdots   \\
s_k & s_{k+1} & \cdots & s_{2k}
\end{array} \right|
\qquad (k = 0,1,\cdots)
$$
are positive.      \\
\indent Assuming such existence (as in [Akh65, (3.4)]), there may or may not be a {\it unique} solution.  With uniqueness, the moment problem is called {\it determinate}; without, it is {\it indeterminate}.  The Stieltjes transform
$$
f(z) := \int_{-\infty}^{\infty} \frac{d \mu(u)}{u - z}
$$
is defined for $z$ in the complex plane cut along the negative real axis ([Akh65, (3.5)]; it was for this that Stieltjes needed to develop his integral, below), and one has the asymptotic expansion
$$
f(z) \sim -(s_0/z + s_1/z^2 + \cdots), \quad (z \to \infty).
$$
In the indeterminate case, there is a parametrization of the solutions due to Nevanlinna in 1922 [Akh65, (3.11)].  The $n$th stage of the continued fraction for $f(z)$ can be expressed via the composition of $n$ M\"obius (bilinear) maps in terms of the `remainder', $f_n(z)$ [Akh, \S 3.3]. \\
\indent When the integration is over the half-line $[0,\infty)$, one has the {\it Stieltjes moment problem}; the term moment problem is due to Stieltjes (1994).  When the integration is over a finite interval ($[0,1]$, say), we have the {\it Hausdorff moment problem} of 1923; it is always determinate.\\
\indent  T. J. Stieltjes (1856-1894) is the leading figure here, with his two-part classic [Stie94/95], the second part posthumous.  He was Dutch, a student at Delft, later an academic at Toulouse.  He had an unconventional career; his most important contact was Hermite, with whom he corresponded extensively; he died sadly young, at 38. \\
\indent Stieltjes's main interest was in continued fractions.  But his work is best remembered now for his introduction of the {\it Stieltjes integral}, which occurs in both Lebesgue-Stieltjes and Riemann-Stieltjes forms, and also for the Stieltjes transform (iterated Laplace transform; Widder [Wid41, Ch. VIII]) and the Stieltjes-Vitali theorem (or Vitali's convergence theorem: [Tit39, \S 5.21]).  Kjeldsen [Kje93] gives a detailed analysis of the evolution of Stieltjes's ideas, referring for the history of the Lebesgue-Stieltjes integral to the Epilogue of Hawkins's book [Haw70].    \\
\indent The reader may have met the Weyl limit-point/limit circle dichotomy for essential self-adjointness (see e.g. [ReeS, 152, 319-320]) in the context of Sturm-Liouville theory.  There is a discrete analogue for Jacobi matrices $\cal{J}$ as in \S 4.5 (cf. Simon [Sim98]).  The limit-{\it circle} case there gives an {\it indeterminate} moment problem [Akh65, Th. 2.1.2]; the limit-{\it point} case gives a {\it determinate} one [Akh65, Cor. 2.2.4] (as cardinality would suggest). \\
\indent {\it Carleman's condition} gives a useful sufficient condition for determinacy: if
$$
\sum_1^{\infty} s_{2n}^{-1/2n} = \infty,
$$
then the (Hamburger) moment problem $(Ham)$ is determinate.  This is the most widely used sufficient condition for determinacy.  It has several equivalent forms, according to the {\it Denjoy-Carleman theorem} (Denjoy, 1921, Carleman 1926 [Car26]).  One involves {\it quasi-analyticity} ([Car26]; Koosis [Koo88, Ch. IV], Rudin [Rud74, Th. 19.11]).  Another involves the {\it Krein condition}: finiteness or otherwise of the {\it logarithmic integral}
$$
\int_{-\infty}^{\infty} \frac{\log Q(x)}{1 + x^2} dx, \qquad 
Q(x) := \sum_0^{\infty} x^n/c_n 
$$
([Rud74]; [Koo88, Ch. IV]).  For this, note that the moment sequence $(c_n)$ satisfies $c_n^2 \leq c_{n-1} c_{n+1}$, by the Cauchy-Schwarz(-Bunyakovskii) inequality.  (This inequality is so close to quasi-analyticity that it can be assumed here without loss of generality, as Rudin explains.) \\
\indent These matters are of great importance in probability theory, in particular in prediction theory (\S 9.6).  For a fine treatment of prediction of Gaussian processes in continuous time, in particular for the work of M. G. Krein in this area (including Krein's theory of strings), see Dym and McKean [DymM76] (cf. Marcus [Marc77]). \\
     
{\bf 6. Mechanical devices} \\

 Throughout his career, Chebyshev was deeply interested in mechanical devices of various sorts, and wrote extensively on them.  We will confine ourselves here to one aspect (below), but note that Chebyshev's breadth of interest and versatility was remarkable even in its time, and would hardly happen in this modern era of increased specialization. \\
\indent The golden age of the steam engine is past, but was enormously influential at the time of the Industrial Revolution, and indeed in Chebyshev's time.  The golden age of the internal combustion engine is now in sight in its turn, but has been enormously influential in our own times.  It is thus salutary to realise that the transference of oscillating linear motion of a piston to circular motion of a wheel is a non-trivial matter.  The device needed is called a {\it linkage} \\
\indent Thomas Newcomen (1664-1729) built his steam engine of 1712 to pump water out of mines.  James Watt (1736-1819) improved this, made it mobile and used it to power a locomotive in 1776.  This gave a powerful acceleration to the Industrial Revolution.  Much of 19th-century history was influenced by steam power, in railways and ships. \\
\indent Watt's linkage was approximate, not perfect (in view of which it is remarkable that his steam engines worked as well as they did!)  Perfect linkages did not emerge for nearly a century.  The {\it Peaucellier-Lipkin linkage} was found independently by Charles-Nicolas Peaucellier (1832-1919) in 1864 and Lipman Israelevich Lipkin (1840-1876) in 1871.  Lipkin was a student of Chebyshev. \\
\indent For two books on linkages, one from the 19th century, one from the 21st, see Kempe [Kem77], O'Rourke [O'Ro11].  See also the numerous papers in this area by the contemporary mathematician S. C. Power. \\

{\bf 7. Pupils and descendants} \\

 Chebyshev had seven pupils (including the brothers A. A. and V. A. Markov), the most distinguished two being A. A. Markov and A. M. Lyapunov.\\
 
{\it Andrei Andreyevich Markov} (1856-1922), PhD 1884 \\
 Markov's PhD thesis was `On certain applications of continued fractions'.  Krein's account of his work (and Chebyshev's) on `limiting values of integrals' has been mentioned in \S 4.8. \\
\indent Markov is of course best remembered for his introduction of Markov chains (1906).  His book [Mark12] of 1912, an early classic, is the German translation of the second (of four) Russian editions.  Seneta [Sen96] gives a good account of the impetus behind the birth of Markov chains.  Until then, different random quantities had often been tacitly taken as `independent by default'.  We take the Markovian condition of `given the present, the past is irrelevant for predicting the future' for granted, but it was a ground-breaking conceptual advance in its time.  He is also remembered for {\it Markov's inequality}, a result similar to Chebyshev's inequality; see e.g. [Bil95, 1.5].  \\
\indent Seneta [Sen84] also gives a good historical account of the central limit theorem and least squares in pre-revolutionary Russia. \\
\indent See [She89] for an assessment of Markov's work on probability, and [Ond81] for an account of the extensive mathematical correspondence between Markov and Chuprov in the period 1910-17. \\
\indent Markov and Sonin (Sonine) edited Chebyshev's Collected Works [MarkS]; for his complete works in Russian, see [CheW]. \\      
\indent Markov had seven pupils, including G. F. Voronoy (PhD 1896), A. S. Besicovitch (PhD 1912), and J. D. Tamarkin (PhD 1917).\\
\indent  Voronoy's two pupils included W. Sierpinski (PhD 1906).  Sierpinski's many pupils included S. Mazurkiewcz, 1913, K. Kuratowski and A. Rajchman, 1921, S. Saks, 1922, A. Zygmund, 1923, and J. Neyman, 1924.  Neyman's many pupils included G. Dantzig and E. Lehmann in 1946 and L. Le Cam in 1952. \\ 
\indent Besicovitch's 9 students included P. A. P. Moran, R. O. Davies and S. J. Taylor (1954), and J. M. Marstrand (1955). \footnote{I met Besicovitch and heard him speak, on Chebyshev, in Cambridge; James Taylor was my first Head of Department, in London.} \\
\indent Tamarkin's 28 students included Nelson Dunford (1936). \\
\indent We note that of this list of 18 of Markov's descendants, the first three were Russian, the next seven Polish; of the last eight, three were US, three UK, one Australian and one French (plus five cases of emigration, to US or UK).  All were good mathematicians; several were great ones: a fine illustration of the quality and world-wide character of Chebyshev's descendants. \\

{\it Alexander Mikhailovich Lyapunov} (1859-1918), PhD 1885 \\
 Lyapunov's thesis was `On the stability of elliptic equilibrium forms of a rotating liquid' (one thinks of the Earth, an oblate spheroid, largely composed of rotating fluid, and wonders about its stability \ldots).  It led him to the work he is best known for, {\it stability of dynamical systems}, in 1892 (the same year as Poincar\'e's work on celestial mechanics and the stability of the solar system).  {\it Lyapunov functions} are all around us; see e.g. Hahn [Hah63], and for probabilistic applications, Menshikov, Popov and Wade [MenPW].  So too are {\it Lyapunov exponents} (characteristic exponents), which arise in products of random matrices; see e.g. Guivarc'h [Gui], Ledrappier [Led], Bougerol and Lacroix [BouL]. \\
\indent Lyapunov proved the central limit theorem rigorously for the first time in 1901, using the method of characteristic functions.  These are indeed characteristic of the Russian school of probability; one need only think of such classics as Ibragimov and Linnik [IbrL] and Petrov [Pet72]. \\
\indent Lyapunov had two pupils, including V. A. Steklov (PhD, Kyiv, 1901), after whom the Steklov Mathematical Institute of the Russian Academy of Sciences is named (on his death in 1926, for his role in founding it in 1921). \\

\indent In addition, M. G. Krein 1907-1989) was a direct descendant, via Korkin, Grave and Chebotarev.  Among his many contributions, his work on the moment problem is touched on in \S 5 and his work on prediction theory in \S 9.6; Krein's theory of strings is developed in [DymM76]; for his work on entire and meromorphic functions see Ostrovskii [Ost94]; the Krein-Milman theorem is one of the cornerstones of functional analysis. \\  
\indent We must add here four illustrious names not in the direct line of descent: N. N. Luzin (a pupil of Egorov), his pupils A. Ya. Khinchin and A. N. Kolmogorov, the last of `the three great names of probability', after Bernoulli and Chebyshev, and S. N. Bernshtein (Bernstein), a pupil of Hilbert. \\

{\bf 8.  Links with the West} \\

 Chebyshev was (for most of his life) a rich man, and like many affluent Russians of his time enjoyed spending time in Europe.  He was also keenly aware that for the full potential impact of his work to be recognised, personal contact with leading Western mathematicians was important, and so was overcoming the language barrier (at least between Russian and French).  Chebyshev was a notoriously bad correspondent, which made personal contact all the more important. \\ 
\indent Chebyshev made a `grand tour' of the West (mainly France and Britain) in 1852, on which he reported to the St. Petersburg Academy.  He visited railways, mines, foundries, mills, and the factories of the manufacturers of large steam engines, for trains and ships. \\
\indent Chebyshev had contacts with Bienaym\'e (\S 3; see [HeyS]).  He had close friendships with Charles Hermite (1822-1901), and E. C. Catalan (1814-1894).  He had good links with Liouville, whom he visited repeatedly in the period 1852-1878, and published 17 papers in Liouville's journal, {\sl J. Math. Pures Appl.}  He published 3 papers in Crelle's journal {\sl J. Reine Angew. Math.}, and 5 in {\sl Acta Mathematica}, after its foundation in 1885. \\
\indent His contacts with the West declined after around 1884.  As Butzer and Jongmans close [ButJ89]: `Maybe Chebyshev had finally found in Russia what he had previously been searching for in Paris, and perhaps Berlin, namely participation in an active mathematical life, now centred around his own students.  The school that he founded was growing steadily and it has since received international fame'. \\  

{\bf 9. Legacy} \\

\indent For reasons of space, we confine ourselves to a few specifics. \\

{\it 1. Generalised (Bienaym\'e-)Chebyshev inequalities} \\
\indent For a variety of generalisations, to multivariate situations, using convexity arguments etc., see e.g. Kingman [Kin63], Whittle [Whi58]. \\

{\it 2. Pad\'e approximation} \\
\indent The Pad\'e approximant of a function $f$ by a rational function -- the ratio of a polynomial of degree $m$ by one of degree $n$ -- is the function whose Taylor series agrees with that of $f$ up to terms of degree $m+n$.  They were developed by Henri Pad\'e (1863-1953) from 1890 on.  They are widely used, in approximation theory and in branches of physics [Bre91]. \\
 
{\it 3. Orthogonal polynomials and stochastic processes} \\
\indent Tridiagonal matrices arise naturally in the theory of birth-and-death processes (for which see e.g. [KarT, Ch. 4]).  Spectral representations are obtained for their transition probabilities in terms of the corresponding orthogonality measure (\S 5) in [KarM].  For an account of further applications of this kind, see Schoutens [Scho00]. \\  

{\it 4.  Markov chains} \\
\indent Markov chains, ubiquitous nowadays, became widely used following their appearance in the first (1950) edition of Feller's book [Fel68], and the books by Kemeny and Snell [KemS] in 1960 and Kemeny, Snell and Knapp [KemSK] in 1966, on the finite and (denumerably) infinite cases.  The more recent theory was given a tremendous boost by the introduction of MCMC (Markov chain Monte Carlo) methods in statistics, between the first (1993) and second (2009) editions of Meyn and Tweedie [MeyT]. \\ 
      
{\it 5. Positive matrices and operators} \\
\indent The theory of positive matrices rests on the Perron-Frobenius theory of O. Perron (1907) and G. Frobenius (1909, 1912).  It has been extensively applied to matrices of transition probabilities of Markov chains.  See Seneta [Sen81] (and the additional bibliography in the 2006 reprinting). \\
\indent The theory can be extended to the continuous setting of the Krein-Rutman theory of positive operators; see e.g. Schaefer [Scha74]. \\
 
{\it 6. Prediction theory} \\
\indent Simon's books [Sim05], [Sim11] give a detailed account of orthogonal polynomials on the unit circle (OPUC -- and on the real line, OPRL, in [Sim11]).  Both OPUC and Szeg\H{o}'s theorem are extremely useful in prediction theory (of stationary stochastic processes); see e.g. the author's survey, [Bin12a].  In discrete time, one works on the unit circle and uses Szeg\H{o}'s theorem; in continuous time, one works on the line and uses the Krein condition (on the logarithmic integral) which corresponds to it; one passes between the unit disc $\mathbb{D}$ and the upper half-plane by the M\"obius map $w = (z - i)/(z + i)$.  For background on the interplay between the two, see [Sim11], the last chapter of Akhiezer [Akh65], the first chapter of Koosis [Koo88], and the end of [Bin12a].\\
\indent In his earlier survey of the moment problem in 1998, Simon remarks [Sim88, p.86] `Non-uniqueness [indeterminacy] only occurs in somewhat pathological situations, but the theory is so elegant and beautiful that it has captivated analysts for a century'.  This is a very good thing for probabilists!  To them the situation is just the reverse: determinacy corresponds in the stochastic-process setting to the future being determined by the remote past, leaving no room for new randomness (`new blood') -- which is degenerate.  The tools developed in analysis have proved extremely useful in probability. \\

{\it 7. Higher dimensions} \\
\indent It is probabilistically of great interest to extend this theory to finite  dimensions [Bin12b], and indeed infinite dimensions [BinM], [Bin20]. \\
\indent For multidimensional continued fractions, see Schweiger [Schw00].  For multidimensional moment problems, see Putinar and Schm\"udgen [PutS]. \\ 

{\bf 10. Conclusion} \\
\indent It is a tribute to the extraordinary breadth, as well as depth, of Chebyshev's work that, in this later epoch of greater specialization, no one can study his work without learning a great deal. \\
\indent In the centre of St Paul's Cathedral in London, built by Sir Christopher Wren after the Great Fire of London of 1666, is the inscription {\it si monumentum requiris, circumspice} (If you seek my monument, look around you).  Chebyshev is the father of Russian mathematics, and of Russian probability in particular.  His monument is all around us. \\ 

\noindent {\bf References} \\
\noindent [Akh56] N. I. Akhiezer, {\sl Theory of approximation}.  Frederick Ungar, 1956; Dover, 1992 (Russian, OGIZ, 1947). \\ 
\noindent [Akh65] N. I. Akhiezer, {\sl The classical moment problem and some related questions in analysis}.  Hafner/Oliver \& Boyd, 1965; SIAM, Classics App. Math. {\bf 82}, 2021 (Russian, 1961). \\
\noindent [Ber45] J. Bertrand, M\'emoire sur le nombre de valeurs que peut prendre une fonction quant on y permute des lettres qu'elle renferme.  {\sl J. Ec. Roy. Polytechnique} {\bf 18} (1845), 123-140. \\
\noindent [Bie53] I. J. Bienaym\'e, Consid\'erations \`a l'appui de la d\'ecouverte de Laplace dans la loi de probabilit\'es dans la m\'ethode des moindres carr\'es.  {\sl Comptes Rendus Acad. Sci. Paris} {\bf 37} (1853), 309-324 (reprinted, {\sl J. Math. Pures Appl.} (2) {\bf 12} (1867), 158-176). \\
\noindent [Bil95] P. Billingsley, {\sl Probability and measure}, 3rd ed., Wiley, 1995 (1st ed. 1979, 2nd ed. 1986). \\
\noindent [Bin12a] N. H. Bingham, Szeg\H{o}'s theorem and its probabilistic descendants.  {\sl Probability Surveys} {\bf 9} (2012), 287-324. \\
\noindent [Bin12b] N. H. Bingham,  Multivariate prediction and matrix Szeg\H{o} theory.  {\sl Probability Surveys} {\bf 9} (2012), 325-339. \\
\noindent [Bin20] N. H. Bingham, Prediction theory for stationary functional time series. arXiv:2011.09937.  \\
\noindent [BinM] N. H. Bingham and Badr Missaoui, Aspects of prediction. {\sl J. Appl. Prob.} {\bf 51A} (2014), 189-201. \\
\noindent [Boa52] R. P. Boas, Review of [Kre51].  {\sl Math. Reviews} {\bf 13} (1952), 445c (MR0044591). \\
\noindent [BouL] P. Bougerol and J. Lacroix, {\sl Products of random matrices with applications to Schr\"odinger operators}.  Birkh\"auser, 1985 ( Progress in Probability and Statistics {\bf 8}). \\   
\noindent [Bre91] C. Brezinski, {\sl History of continued fractions and Pad\'e approximants}.  Springer, 1991. \\
\noindent [ButJ89] P. Butzer and F. Jongmans, P. L. Chebyshev (1821-1894) and his contacts with Western European scientists.  {\sl Historia Mathematica} {\bf 16} (1989), 46-68. \\
\noindent [ButJ99] P. Butzer and F. Jongmans, P. L. Chebyshev (1821-1894): A guide to his life and works.  {\sl J. Approximation Theory} {\bf 96} (1999), 111-138. \\
\noindent [Car26] T. Carleman, {\sl Les fonctions quasi-analytiques}.  Gauthier-Villars, 1926. \\
\noindent [Che46] P. L. Chebyshev, D\'emonstration \'el\'ementaire d'une proposition g\'en\'erale de la th\'eorie des probabilit\'es.  { \sl J. Reine Angew. Math.} {\bf 33} (1846), 259-267 (Works I, 15-26). \\
\noindent [Che52] P. L. Chebyshev, M\'emoire sur les nombres premiers.  {\sl J. Math. Pures Appl.} {\bf 17} (1852), 366-390. \\
\noindent [Che58] P. L. Chebyshev, Sur les fractions continues.  {\sl J. Math.} (2) {\bf 3} (1858), 289-323 (Oeuvres 1, 201-230).  Read on 12.1.1853. \\
\noindent [Che59] P. L. Chebyshev, Sur les questions de minima qui se rattechent \`a la r\'epresentation approximative des fonctions.  {\sl M\'em. Acad. St. P\'etersb.} {\bf 7}(6) (1859), 199-291.  Read on 9.10.1857. \\
\noindent [Che67] P. L. Chebyshev, Des valeurs moyennes.  {\sl J. Math. Pures Appl.} (2) {\bf 12} (1867), 177-184 (Works I, 685-694). \\
\noindent [Che74] P. L. Chebyshev, Sur les valeurs limites des int\'egrales.  {\sl J. Math. Pures Appl.} (2) {\bf 19} (1874), 157-160. \\
\noindent [Che87] P. L. Chebyshev, Sur deux th\'eor\`emes relatifs aux probabilit\'es.  {\sl Acta Math.} {\bf 14} (1890-1891), 305-315 (Russian, 1887; Works II, 479-491). \\
\noindent [CheW] {\sl Complete Collected Works} (Russian, 1946/1951).  Izdatelstvo Akad. Nauk SSR, Moscow/Leningrad, Vol. I (1946), 342p, II (1947), 520p, III (1948), 414p, IV (1948), 255p, V (1951), 474p. \\
\noindent [Dav63] P. J. Davis, {\sl Interpolation and approximation}.  Blaisdell, 1963; Dover, 1975. \\
\noindent [DymM76] H. Dym and H. P. McKean, {\sl Gaussian processes, function theory, and the inverse spectral problem}.  Academic Press, 1976. \\
\noindent [Fel68] W. Feller, {\sl An introduction to probability theory and its applications}, Vol. I, 3rd ed., Wiley, 1968 (1st ed. 1950, 2nd ed. 1957). \\
\noindent [GneK] B. V. Gnedenko and A. N. Kolmogorov, {\sl Limit distributions for sums of independent random variables} (translated by K.-L. Chung, with an Appendix by J. L. Doob).  Addison-Wesley, 1954 (Russian, 1949). \\
\noindent [Gui] Y. Guivarc'h, Quelques propri\'et\'es asymptotiques des produits de matrices al\'eatoires.  {\sl Ecole d'Et\'e de Probabilit\'es de Saint-Flour VIII-1978}, 177-250 ({\sl Lecture Notes in Math.} {\bf 774}, Springer, 1980). \\
\noindent [Hah63] W. Hahn, {\sl Theory and application of Liapunov's direct method}.  Prentice-Hall, 1963. \\
\noindent [Har05] G. H. Hardy, {\sl The integration of functions of a single variable}.  Cambridge University Press, 1905 (2nd ed. 1916;  Cambridge Tracts Math. {\bf 2}).\\
\noindent [HarW] G. H. Hardy and E. M. Wright, {\sl An introduction to the theory of numbers}, 6th ed. (rev. D. R. Heath-Brown and J. H. Silverman), Oxford University Press, 2008. \\
\noindent [Haw70] T. Hawkins, {\sl Lebesgue's theory of integration: Its origin and development}.  University of Wisconsin Press, 1970. \\  
\noindent [HeyS] C. C. Heyde and E. Seneta, {\sl I. J. Bienaym\'e: Statistical theory anticipated}.  Springer, 1977. \\
\noindent [IbrL] I. A. Ibragimov and Yu. V. Linnik, {\sl Independent and stationary sequences of random variables}.  Wolters-Noordhoff, 1971. \\   
\noindent [Kah85] J.-P. Kahane, {\sl Some random series of functions}, 2nd ed.  Cambridge University Press, 1985 (1st ed., D. C. Heath, 1968). \\
\noindent [Kar66] S. Karlin, {\sl Tchebycheff systems: with applications in analysis and statistics}.  Interscience, 1966. \\
\noindent [KarM] S. Karlin and J. L. McGregor, The differential equations of birth and death processes and the Stieltjes moment problem.  {\sl Trans. Amer. Math. Soc.} {\bf 85} (1957), 486-546. \\
\noindent [KarT] S. Karlin and H. M. Taylor, {\sl A first course in stochastic processes}, 2nd ed., Academic Press, 1975. \\
\noindent [Kem77] A. B. Kempe, {\sl How to draw a straight line: A lecture on linkages}.  Macmillan, 1877. \\
\noindent [KemS] J. G. Kemeny and J. L. Snell, {\sl Finite Markov chains}.  Van Nostrand, 1960. \\
\noindent [KemSK] J. G. Kemeny, J. L. Snell and A. W. Knapp, {\sl Denumerable Markov chains}.  Van Nostrand, 1966. \\
\noindent [Khr08] S. Khrushchev, {\sl Orthogonal polynomials and continued fractions, from Euler's point of view}.  Cambridge University Press, 2008 (Encycl. Math. Appl. {\bf 122}). \\
\noindent [Kin63] J. F. C. Kingman, On inequalities of the Tchebychev type.  {\sl Proc. Cambridge Phil. Soc.} {\bf 59} (1963), 135-146. \\
\noindent [Kje93] T. H. Kjeldsen, The early history of the moment problem.  {\sl Historia Mathematica} {\bf 20} (1993), 19-44. \\ 
\noindent [Koo88] P. Koosis, {\sl The logarithmic integral I}.  Cambridge University Press, 1988. \\
\noindent [Kre51] M. G. Krein, Chebyshev's and Markov's ideas in the theory of limiting value of integrals and their further development (in Russian).  {\sl Usp. Mat. Nauk} {\bf 44} (1951), 3-120. \\
\noindent [Kry36] A. N. Krylov (ed.), Theory of probability (Russian).  Lectures 1879-80 by P. L. Chebyshev, from notes taken by A. M. Lyapunov.  Nauka, Moscow/Leningrad, 1936. \\
\noindent [Lan09] E. Landau, {\sl Handbuch der Lehre von der Verteilung der Primzahlen}, 2nd ed., Chelsea, 1953 (1st ed. 1909). \\
\noindent [Led] F. Ledrappier, Quelques propri\'et\'es des exposants characteristiques.   {\sl Ecole d'Et\'e de Probabilit\'es de Saint-Flour XII-1982}, 305-396 ({\sl Lecture Notes in Math.} {\bf 1097}, Springer, 1984). \\  
\noindent [Lor53] G. G. Lorentz, {\sl Bernstein polynomials}.  University of Toronto Press, 1953 (Math. Expos. {\bf 8}). \\ 
\noindent [Lor86] G. G. Lorentz, {\sl Approximation of functions}, 2nd ed.  Chelsea, 1986 (1st ed. Holt, Rinehart and Winston, 1966). \\
\noindent [Mai74] L. E. Maistrov, {\sl Probability theory: A historical sketch}.  Academic Press, 1974. \\
\noindent [Marc77] M. B. Marcus, Book review, [DymM76].  {\sl Ann. Probab.} {\bf 5} (1977), 1047-1051. \\
\noindent [Mark12] A. A. Markov, {\sl Wahrscheinlichkeitsrechnung}.  B. G. Teubner, 1912. \\
\noindent [MarkS] A. Markoff and N. Sonin, {\sl Oeuvres de P. L. Tchebyshef}, 2 volumes, St. Petersburg, 1899/1907 (repr. Chelsea, New York, 1952). \\
\noindent [MenPW] M. Menshikov, S. Popov and A. Wade, {\sl Non-homogeneous random walks: Lyapunov function methods for near-critical stochastic systems}.  Cambridge University Press, 2017 ( Cambridge Tracts Math. {\bf 209}). \\
\noindent [MeyT] S. Meyn and R. L. Tweedie, {\sl Markov chains and stochastic stability}, 2nd ed.  Cambridge University Press, 2009 (1st ed. Springer, 1993). \\
\noindent [Nat64] I. P. Natanson, {\sl Constructive function theory, I: Uniform approximation}.  Frederick Ungar, 1964. \\
\noindent [Nat65a] I. P. Natanson, {\sl Constructive function theory, II: Approximation in mean}.  Frederick Ungar, 1965. \\
\noindent [Nat65b] I. P. Natanson, {\sl Constructive function theory, III: Interpolation and approximation quadratures}.  Frederick Ungar, 1965. \\
\noindent [Ond81] Kh. H. Ondar (ed.), {\sl The correspondence between A. A. Markov and A. A. Chuprov on the theory of probability and mathematical statistics}.  Springer, 1981. \\
\noindent [O'Ro11] J. O'Rourke, {\sl How to fold it: The mathematics of linkages, origami and polyhedra}.  Cambridge University Press, 2011. \\
\noindent [Ost94] I. V. Ostrovskii, The investigations of M. G. Krein in the theory of entire and meromorphic functions and their further development.  {\sl Ukr. Mat. Zh.} {\bf 46} (1994), 87-100. \\ 
\noindent [Pet72] V. V. Petrov, {\sl Sums of independent random variables}.  Springer, 1972 
(Ergebnisse Math. {\bf 82}). \\
\noindent [PutS] M. Putinar and K. Schm\"udgen, Multidimensional determinateness.  {\sl Indiana U. Math. J.} {\bf 57} (2008), 2931-2968. \\
\noindent [Ram19] S. Ramanujan, A proof of Bertrand's postulate.  {\sl J. Indian Math. Soc.} {\bf 11} (1919), 181-182 ({\sl Collected Papers of Srinivasa Ramanujan} 208-209, Cambridge University Press, 1927). \\
\noindent [ReeS] M. Reed and B. Simon, {\sl Methods of modern mathematical physics, II: Fourier analysis, self-adjointness}.  Academic Press, 1975. \\
\noindent [Rit48] J. F. Ritt, {\sl Integration in finite terms: Liouville's theory of elementary methods}.  Columbia University Press, 1948. \\
\noindent [RogJ] C. A. Rogers and J. E. Jayne, $K$-analytic sets.  Part I (p. 1-181) in {\sl Analytic sets} (C. A. Rogers et al.), Academic Press, 1980.\\
\noindent [Ros88] H. E. Rose, {\sl A course in number theory}. Oxford University Press, 1988. \\
\noindent [Rud74] W. Rudin, {\sl Real and complex analysis}, 2nd ed., McGraw-Hill, 1974 (1st ed. 1966). \\
\noindent [Scha74] H. H. Schaefer, {\sl Banach lattices and positive operators}.  Springer, 1974 (Grundl. math. Wiss. {\bf 215}). \\
\noindent [Scho00] W. Schoutens, {\sl Stochastic processes and orthogonal polynomials}.  Springer, 2000 (Lecture Notes in Statistics {\bf 146}). \\
\noindent [Schw00] F. Schweiger, {\sl Multidimensional continued fractions}.  Oxford University Press, 2000. \\
\noindent [Sea67] H. L. Seal, Studies in the history of probability and statistics XV: The historical development of the Gauss linear model. {\sl Biometrika} {\bf 54} (1967), 1-24. \\
\noindent [Sen81] E. Seneta, {\sl Non-negative matrices and Markov chains}.  Springer, 1981.\\
\noindent [Sen84] E. Seneta, The central limit theorem and linear least squares in pre-revolutionary Russia: The background.  {\sl Math. Scientist} {\bf 9} (1984), 37-77. \\
\noindent [Sen96] E. Seneta, Markov and the birth of chain dependence theory.  {\sl Internat. Stat. Review} {\bf 64} (1996), 255-263. \\
\noindent [She89] O. B. Sheynin, A. A. Markov's work on probability.  {\sl Arch. Hist. Exact Sciences} {\bf 39} (1989), 337-377. \\
\noindent [She94] O. B. Sheynin, Chebyshev's lectures on the theory of probability.  {\sl Arch. Hist. Exact Sciences} {\bf 46} (1994), 321-340. \\
\noindent [ShoT] J. A. Shohat and J. D. Tamarkin, {\sl The problem of moments}.  Amer. Math. Soc., 1943 (Math. Surveys {\bf 1}). \\
\noindent [Sie28]. W. Sierpinski, {\sl Le\c cons sur les nombres transfinis}.  Gauthier-Villars, 1928 (Collection Borel). \\
\noindent [Sie64] W. Sierpinski. {\sl Elementary theory of numbers}.  PWN, 1964.\\
\noindent [Sim98] B. Simon, The classical moment problem as a self-adjoint finite-difference operator.  {\sl Adv. Math.} {\bf 137} (1998), 82-203. \\
\noindent [Sim05]  B. Simon, {\sl Orthogonal polynomials on the unit circle.  Part I, Classical theory; Part II, Spectral theory}.  Amer. Math. Soc., 2005 (AMS Colloq. Publ. {\bf 54} Parts 1,2). \\
\noindent [Sim11] B. Simon, {\sl Szeg\H{o}'s theorem and its descendants}.  Princeton University Press, 2011. \\
\noindent [Stie94/95]  Recherches sur les fractions continues I, II.  {\sl Ann. Fac. Sci. U. Toulouse} {\bf 8} (1894), 1-122, {\bf 9} (1895), 5-47 (available online).  \\
\noindent [Stig86] S. M. Stigler, {\sl The history of statistics: The measurement of uncertainty before 1900}.  Harvard University Press, 1986. \\
\noindent [Sze59] G. Szeg\H{o}, {\sl Orthogonal polynomials}.  Amer. Math. Soc., 1959 (AMS Colloq. Publ. {\bf XXIII}). \\
\noindent [Ten15] G. Tenenbaum, {\sl Introduction to analytic and probabilistic number theory}, 3rd ed.  Amer. Math. Soc., 2015 (Grad. Studies Math. {\bf 163}). \\
\noindent [Tit39] E. C. Titchmarsh, {\sl The theory of functions}, 2nd ed., Oxford University Press, 1939 (1st ed. 1932). \\
\noindent [Whi58] P. Whittle, A multivariate generalization of Tchebychev's inequality.  {\sl Quart. J. Math.} {\bf 9} (1958), 232-240. \\
\noindent [Wid41] D. V. Widder, {\sl The Laplace transform}.  Princeton University Press, 1941. \\

\noindent N. H. Bingham, Mathematics Department, Imperial College, London SW7 2AZ, UK; n.bingham@ic.ac.uk \\ 

\end{document}